\documentclass[11pt,a4paper]{article}
\usepackage[cp1251]{inputenc}
\usepackage{amsfonts, amssymb, amsthm, amsmath}
\usepackage{mathtext}
\usepackage{longtable}
\usepackage[russian]{babel}

\theoremstyle{remark}

\renewcommand{\r}{\mathbb R}
\newcommand{\R}{\mathbb R}

\renewcommand{\le}{\leqslant}
\renewcommand{\ge}{\geqslant}

\renewcommand{\phi}{\varphi}

\newcommand{\eqd}{\stackrel{d}{=}}
\newcommand{\sint}{\text{\footnotesize ${\displaystyle\int}$}}

\title{On asymmetric generalization of the Weibull distribution by
scale-location mixing of normal laws\thanks{Research supported by
the Russian Foundation for Basic Research (project 15-07-04040).}}
\author{Victor Korolev\thanks{Faculty of Computational Mathematics
and Cybernetics, Lomonosov Moscow State University; Institute of
Informatics Problems, Federal Research Center <<Computer Science and
Control>> of the Russian Academy of Sciences; vkorolev@cs.msu.ru},
Lily Kurmangazieva\thanks{Faculty of Computational Mathematics and
Cybernetics, Lomonosov Moscow State University;
lkurmangazieva@gmail.com}, Alexander Zeifman\thanks{Vologda State
University; Institute of Informatics Problems, Federal Research
Center <<Computer Science and Control>> of the Russian Academy of
Sciences; Institute of Socio-Economic Development of Territories,
Russian Academy of Sciences; a$\_$zeifman@mail.ru}}
\date{}

\begin{document}

\maketitle

{\small

{\bf Abstract:} Two approaches are suggested to the definition of
asymmetric generalized Weibull distribution. These approaches are
based on the representation of the two-sided Weibull distributions
as variance-mean normal mixtures or more general scale-location
mixtures of the normal laws. Since both of these mixtures can be
limit laws in limit theorems for random sums of independent random
variables, these approaches can provide additional arguments in
favor of asymmetric two-sided Webull-type models of statistical
regularities observed in some problems related to stopped random
walks, in particular, in problems of modeling the evolution of
financial markets.

\smallskip

{\bf Key words:} Webull distribution; scale-location mixture of
normal laws; variance-mean normal mixture; random sum; stable
distribution; exponential distribution

}

\section{Introduction. The Weibull distribution}

In probability theory and mathematical statistics it is conventional
to use the term {\it Weibull distribution} for a special absolutely
continuous probability distribution concentrated on the nonnegative
half-line with exponential-power type decrease of the tail. It is
called so after the Swedish scientist Waloddi Weibull (1887--1979)
who suggested in 1939 to use this distribution in the study of the
strength of materials \cite{Weibull1939a, Weibull1939b} and
thoroughly analyzed this distribution in 1951 \cite{Weibull1951}
demonstrating good perspectives of the application of this
distribution to the description of many observed statistical
regularities.

Let $\gamma>0$. The distribution of the random variable
$W_{\gamma}$:
$$
{\sf
P}\big(W_{\gamma}<x\big)=\big[1-e^{-x^{\gamma}}\big]\mathbf{1}(x\ge
0),\ \ \ x\in\mathbb{R},\eqno(1)
$$
is called the {\it Weibull distribution} with shape parameter
$\gamma$ (here and in what follows the symbol $\mathbf{1}(C)$
denotes the indicator function of a set $C$).

However, Weibull was not the first to introduce this distribution.
This distribution was for the first time described in 1927 by
Maurice Fr{\'e}chet \cite{Frechet1927} within the context of the
study of the asymptotic behavior of extreme order statistics.
Sometimes distribution (1) is called the {\it Rosin--Rammler
distribution} after Paul Rosin and Erich Rammler, German scientists
who were the first to use this distribution as a model of
statistical regularities in the coal particles sizes in 1933
\cite{RosinRammler1933}. However, this term also does not completely
correspond to the historical truth. In the paper \cite{Stoyan2013}
Dietrich Stoyan directly writes that this distribution was found by
Paul Rosin, Erich Rammler, Karl Sperling \cite{RosinRammler1933,
RosinRammlerSperling1933} and John Godolphin Bennett
\cite{Bennett1936} within the context of particle size. It is well
known that the Weibull distribution family is closed with respect to
the operation of taking minimum of independent random variables. As
it was demonstrated in \cite{Frechet1927, FisherTippett1928}, due to
this property the family of Weibull distributions is one of possible
limit laws for extreme order statistics. B. V. Gnedenko found
necessary and sufficient conditions for the convergence of the
distributions of extreme order statistics to the Weibull
distribution under linear normalization \cite{Gnedenko1943}.
Therefore, this distribution is sometimes called the {\it
Weibull--Gnedenko distribution} \cite{JohnsonKotzBalakrishnan1994}.

The case of small values of the parameter $\gamma\in(0,1]$ is of
special interest for financial and some other applications, since
Weibull distributions with such parameters (sometimes called {\it
stretched exponential distributions} \cite{LaherrereSornette1998,
Sornette_et_al2005, Sornette_et_al2006}) occupy an intermediate
position between distributions with exponentially decreasing tails
(such as exponential and gamma-distributions) and heavy-tailed
Zipf--Pareto-type distributions with power-type decrease of tails.

In the paper \cite{Stoyan2013} cited above, D. Stoyan notes that the
name {\it exponential power distribution} would better fit to
distribution (1), however, the latter term has been occupied by
another absolutely continuous distribution with a similar behavior
of tails \cite{BoxTiao1973, Korolevetal2012, GrigoryevaKorolev2013},
which, unlike distribution (1), has the exponential power {\it
density}
$$
\ell_{\gamma}(x)=\frac{\gamma}{2\Gamma({\textstyle\frac{1}{\gamma}})}e^{-|x|^{\gamma}},\
\ \ x\in\mathbb{R},
$$
with $\gamma>0$, whereas distribution (1) has the exponential power
{\it cumulative distribution function}.

In this paper, for definiteness, for distribution (1) we will use
the traditional term {\it Weibull distribution}.

The Weibull distribution is a special generalized
gamma-distribution. It is widely used in survival analysis
\cite{JohnsonJohnson1999}, in life insurance as a model of the
lifetime distribution, in risk insurance as a model of the claim
size distribution \cite{HoggKlugman1983}, in economics and financial
mathematics as a model of asset returns distribution
\cite{DAddario1974, MittnikRachev1989, MittnikRachev1993} and income
distribution \cite{Bartels1977, Bordleyetal1996}, in reliability
theory as a model of the distribution of time between failures
\cite{Abernethy2004, Lawless1982}, in industrial technology as a
model of the distribution of duration of technological stages or
time intervals between technological changes \cite{SharifIslam1980},
in coal industry for the description of statistical regularities of
particle sizes \cite{RosinRammler1933}, in radio engineering and
radiolocation, in meteorology, hydrology and many other fields, see,
e. g., \cite{JohnsonKotzBalakrishnan1994, Abernethy2004,
Lawless1982, KotzNadaraja2000, JohnsonKotz1970}.

In particular, in \cite{MittnikRachev1989} it was discovered that
this distribution provides the best fit among others to the observed
statistical regularities of the index S\&P500, if its positive and
negative increments are considered separately thus leading to the
concept of a two-sided Weibull distribution with both positive and
negative tails decreasing as an exponential power function. Some
authors suggested to use the Weibull distribution as the errors
distribution in range data modelling \cite{Chenetal2008} or the
distribution of trading duration \cite{EngleRussell1998}.

In \cite{Sornette2000} it was proposed to use a symmetric two-sided
Weibull distribution as an unconditional return distribution. A
symmetric two-sided Weibull distribution was also mentioned in
\cite{MalevergneSornette2004}, but its properties were not explored.
It should be noted that in these papers the attempts to introduce
the two-sided Weibull distribution were rather formal and
descriptive. In these papers as well as in \cite{ChenGerlach2011,
ChenGerlach2013} the elementary properties of these models were
described.

In the present paper, two new approaches are proposed to the
definition of general asymmetric two-sided Weibull distribution by
representing them as variance-mean and more general scale-location
mixtures of normal laws.

Among all scale-location mixtures of normal laws, normal
variance-mean mixtures occupy a special position. A distribution
function $F(x)$ is called a {\it normal variance-mean mixture}, if
it has the form
$$
F(x)=\sint_{\!\!0}^{\infty}\Phi\Big(\frac{x-\beta-\alpha
z}{\sigma\sqrt{z}}\Big)\,dG(z),\ \ \ x\in\r.
$$
with some $\beta\in\mathbb{R}$, $\alpha\in\mathbb{R}$,
$\sigma\in(0,\infty)$, where $\Phi(x)$ is the standard normal
distribution function, $G(z)$ is a distribution function such that
$G(0)=0$. The class of normal variance-mean mixtures is very wide
and contains, say, generalized hyperbolic laws \cite{BN1977, BN1978,
BN1982} and generalized variance gamma distributions
\cite{KorolevZaks2013} which proved to provide excellent fit to
statistical data in various fields from atmospheric turbulence to
financial markets.

In normal variance-mean mixtures, mixing is performed with respect
to both location and scale parameters. But since these parameters
are tightly linked so that the location parameters ({\it means}) are
proportional to the {\it variances} of the mixed laws, actually this
is a one-parameter mixture.

Normal mixtures are limit distributions for sums of a random number
of random variables. Therefore, these approaches can provide
additional grounds for the adequacy of asymmetric two-sided Weibull
distributions in practical problems modeled by stopped random walks,
in particular, related to the description of the evolution of
financial indexes.

The paper is organized as follows. Section 2 contains some auxiliary
results dealing with product representations for Weibull-distributed
random variables by normally and exponentially distributed random
variables. In section 3 similar results are obtained for symmetric
two-sided Weibull distributions. In section 4 the representation is
obtained for the formal asymmetric two-sided Weibull distribution as
a scale-location mixture of normal laws. Finally, in section 5 a new
generalization of the Weibull distribution by variance-mean mixing
of normal laws is proposed. Here a limit theorem for random sums of
independent random variables is also proved illustrating that such
Weibull-type distribution can be reasonable asymptotic
approximation.

\section{Product representations for Weibull-distributed random
variables by normally and exponentially distributed random
variables}

Here we will show that the Weibull distribution admits several
representations in the form of scale mixtures of some well-known
probability distributions. For this purpose here some product
representations for Weibull-distributed random variables by normally
and exponentially distributed random variables will be proved. These
representations will be used in the next section to prove that for
$\gamma\in(0,1]$ the Weibull distribution is a mixed half-normal
distribution and hence, it can be the limiting law for maximum
random sums with finite variances.

Most results presented below actually concern special mixture
representations for probability distributions. However, without any
loss of generality, for the sake of visuality and compactness of
formulations and proofs we will formulate the results in terms of
the corresponding random variables assuming that all the random
variables mentioned in what follows are defined on the same
probability space $(\Omega,\,\mathfrak{A},\,{\sf P})$.

It is obvious that $W_1$ is the random variable with the standard
exponential distribution: ${\sf P}(W_1<x)=\big[1-e^{-x}\big]{\bf
1}(x\ge0)$. The Weibull distribution with $\gamma=2$, that is, ${\sf
P}(W_2<x)=\big[1-e^{-x^2}\big]{\bf 1}(x\ge0)$ is called the Rayleigh
distribution \cite{Rayleigh1880}.

The random variable with the standard normal distribution function
$\Phi(x)$ will be denoted $X$,
$$
{\sf
P}(X<x)=\Phi(x)=\frac{1}{\sqrt{2\pi}}\sint_{\!\!-\infty}^{x}e^{-z^2/2}dz,\
\ \ \ x\in\mathbb{R}.
$$
Let $\Psi(x)$, $x\in\mathbb{R}$, be the distribution function of the
maximum of the standard Wiener process on the unit interval,
$\Psi(x)=2\Phi\big(\max\{0,x\}\big)-1$, $x\in\mathbb{R}$. It is easy
to see that $\Psi(x)={\sf P}(|X|<x)$. Therefore, sometimes $\Psi(x)$
is said to determine the {\it half-normal} distribution.

The symbol $\eqd$ denotes the coincidence of distributions.

\smallskip

{\sc Lemma 1.} {\it The relation}
$$
W_1\eqd\sqrt{2W_1} |X|\eqno(2)
$$
{\it holds, where the random variables on the right-hand side are
independent.}

\smallskip

{\sc Proof}. For $x>0$ we have
$$
{\sf P}\big(|X|\sqrt{W_1}<x\big)={\sf
E}\Psi\big(x/\sqrt{W_1}\big)=2{\sf E}\Phi\big(x/\sqrt{W_1}\big)-1=
2\sint_{\!\!0}^{\infty}\Phi\big(x/\sqrt{z}\big)d[1-e^{-z}]-1=
$$
$$
=2\sint_{\!\!0}^{\infty}\Big[\frac{1}{2}+
\frac{1}{\sqrt{2\pi}}\sint_{\!\!0}^{x/\sqrt{z}}e^{-u^2/2}du\Big]e^{-z}dz-1=
\frac{\sqrt{2}}{\sqrt{\pi}}\sint_{\!\!0}^{\infty}\sint_{\!\!0}^{x/\sqrt{z}}e^{-u^2/2-z}dudz=
$$
$$=
\frac{\sqrt{2}}{\sqrt{\pi}}\sint_{\!\!0}^{\infty}\sint_{\!\!0}^{x^2/u^2}e^{-z}dz
e^{-u^2/2}du=
\frac{\sqrt{2}}{\sqrt{\pi}}\sint_{\!\!0}^{\infty}\Bigl(
1-\exp\Bigl\{-\frac{x^2}{u^2}\Bigr\}\Bigr)e^{-u^2/2}du=
$$
$$=
1-\frac{\sqrt{2}}{\sqrt{\pi}}\sint_{\!\!0}^{\infty}
\exp\Bigl\{-\frac{u^2}{2}-\frac{x^2}{u^2}\Bigr\}du=
1-e^{-\sqrt{2}x}={\sf P}(W_1<\sqrt{2}x),
$$
see, e. g., \cite{GradsteinRyzhik1971}, formula 3.325. This is
nothing else that the exponential distribution with parameter
$\sqrt{2}$. The lemma is proved.

\smallskip

It is easy to see that if $\gamma>0$ and $\gamma'>0$, then ${\sf
P}(W_{\gamma'}^{1/\gamma}\ge x)={\sf P}(W_{\gamma'}\ge
x^{\gamma})=e^{-x^{\gamma\gamma'}}={\sf P}(W_{\gamma\gamma'}\ge x)$,
$x\ge 0$, that is, for any $\gamma>0$ and $\gamma'>0$
$$
W_{\gamma\gamma'}\eqd W_{\gamma'}^{1/\gamma}.\eqno(3) 
$$
In particular, for any $\gamma>0$ we have
$$
W_{\gamma}\eqd W_1^{1/\gamma}.\eqno(4)
$$

\section{Product representations for Weibull-distributed random
variables by random variables with stable distributions}

The distribution function and density of the strictly stable
distribution with the characteristic exponent $\alpha$ and parameter
$\theta$, corresponding to the characteristic function
$$
\mathfrak{f}_{\alpha,\theta}(t)=\exp\big\{-|t|^{\alpha}\exp\{-{\textstyle\frac12}i\pi\theta\alpha\,\mathrm{sign}t\}\big\},\
\ \ \ t\in\r,\eqno(5)
$$
with $0<\alpha\le2$ and $|\theta|\le\min\{1,\frac{2}{\alpha}-1\}$
will be respectively denoted $G_{\alpha,\theta}(x)$ and
$g_{\alpha,\theta}(x)$, $x\in\mathbb{R}$ (see, e. g.,
\cite{Zolotarev1983}).

From (5) it follows that the characteristic function of a symmetric
($\theta=0$) strictly stable distribution has the form
$$
\mathfrak{f}_{\alpha,0}(t)=e^{-|t|^{\alpha}},\ \ \ t\in\r. \eqno(6)
$$

\smallskip

{\sc Lemma 2.} {\it Any symmetric strictly stable distribution with
the characteristic exponent $\alpha$ is a scale mixture of normal
laws with the one-sided strictly stable law $(\theta=1)$ with the
characteristic exponent $\alpha/2$ as the mixing distribution}:
$$
G_{\alpha,0}(x)=\sint_{\!\!0}^{\infty}\Phi\big(x/\sqrt{z}\big)dG_{\alpha/2,1}(z),\
\ \ x\in\r.\eqno(7)
$$

\smallskip

{\sc Proof}. See, e. g., theorem 3.3.1 in \cite{Zolotarev1983}.

\smallskip

In order to prove that any Weibull distribution with
$\gamma\in(0,1]$ is a scale mixture of half-normal laws, we first
make sure that any Weibull distribution with $\gamma\in(0,2]$ is a
scale mixture of Rayleigh distributions.

Everywhere in what follows for $\gamma\in(0,1)$ we will use the
notation $V_{\gamma}=2S_{\gamma,1}^{-1}$, where $S_{\gamma,1}$ is a
random variable with one-sided strictly stable density
$g_{\gamma,1}(x)$.

\smallskip

{\sc Lemma 3.} {\it For any $\gamma\in(0,2]$ the product
representation}
$$
W_{\gamma}\eqd W_2 \sqrt{V_{\gamma/2}}\eqno(8)
$$
{\it holds, where the random variables on the right-hand side of
$(8)$ are independent.}

\smallskip

{\sc Proof}. Write relation (7) in terms of characteristic functions
with the account of (6):
$$
e^{-|t|^{\alpha}}=\sint_{\!\!0}^{\infty}\exp\{-{\textstyle\frac12}t^2z\}g_{\alpha/2,1}(z)dz,
\ \ \ t\in\mathbb{R}.\eqno(9)
$$
Formally setting in (9) $|t|=x$ with $x\ge0$ being an {\it arbitrary
nonnegative} number, we obtain
$$
{\sf
P}(W_{\gamma}>x)=e^{-x^{\gamma}}=\sint_{\!\!0}^{\infty}\exp\{-{\textstyle\frac12}x^2z\}g_{\gamma/2,1}(z)dz.\eqno(10)
$$
At the same time, we obviously have
$$
{\sf P}\big(W_2 \sqrt{V_{\gamma/2}}>x\big)={\sf
P}\big(W_2>x\sqrt{{\textstyle\frac12}S_{\gamma/2,1}}\big)=
\sint_{\!\!0}^{\infty}\exp\{-{\textstyle\frac12}x^2z\}g_{\gamma/2,1}(z)dz.\eqno(11)
$$
Since the right-hand sides of (10) and (11) identically (in $x\ge0$)
coincide, the left-hand sides of these relations identically
coincide as well. The lemma is proved.

\smallskip

{\sc Lemma 4.} {\it For any $\gamma\in(0,1]$ the Weibull
distribution with parameter $\gamma$ is a mixed exponential
distribution}:
$$
{\sf P}(W_{\gamma}>x)=
\sint_{\!\!0}^{\infty}e^{-\frac12 zx}g_{\gamma,1}(z)dz,\ \ \ x\ge0.
$$

\smallskip

{\sc Proof.} From $(3)$ it follows that $W_2\eqd\sqrt{W_1}$.
Therefore lemma 3 implies that for $\gamma\in(0,1]$
$$
W_{2\gamma}\eqd W_2 \sqrt{V_{\gamma}}\eqd\sqrt{W_1 V_{\gamma}}
$$
or, with the account of $(3)$,
$$
W_{\gamma}\eqd W_{2\gamma}^2\eqd W_1 V_{\gamma},\eqno(12)
$$
that is,
$$
e^{-x^{\gamma}}={\sf P}(W_{\gamma}>x)={\sf
P}(W_1>{\textstyle\frac12}S_{\gamma,1}x)=
\sint_{\!\!0}^{\infty}e^{-\frac12 zx}g_{\gamma,1}(z)dz,\ \ \ x\ge0.
$$
The lemma is proved.

\smallskip

{\sc Corollary 1.} {\it Any Weibull distribution with
$\gamma\in(0,1]$ is infinitely divisible.}

\smallskip

{\sc Proof.} This statement is actually a particular case of a
result on generalized gamma-convolutions in \cite{Bondesson1979}
(see corollary 2 there). However, we can also see that the desired
assertion immediately follows from (12) and the result of Goldie
\cite{Goldie1967} stating that the product of two independent
non-negative random variables is infinitely divisible if one of the
two is exponentially distributed.

\smallskip

{\sc Theorem 1}. {\it For any $\gamma\in(0,1]$ the Weibull
distribution with parameter $\gamma$ is a scale mixture of
half-normal laws}:
$$
W_{\gamma}\eqd |X| \sqrt{2W_1V_{\gamma}^2},\eqno(13)
$$
{\it where the random variables on the right-hand side of $(13)$ are
independent.}

\smallskip

{\sc Proof.} To prove the desired result it suffices to use the
representation for $W_1$ proved in lemma 1 on the right-hand side of
(12) and obtain (13). The theorem is proved.

\smallskip

Theorem 1 makes it possible to represent the Weibull distribution
with $\gamma\in(0,1]$ as
$$
{\sf P}(W_{\gamma}<x)={\sf
E}\Psi\big(x/\sqrt{2W_1V_{\gamma}^2}\big)=\sint_{\!\!0}^{\infty}\Psi\big(x/\sqrt{y}\big)dH_{\gamma}(y),\
\ \ x\in\mathbb{R},
$$
where
$$
H_{\gamma}(y)={\sf P}(2W_1V_{\gamma}^2<y)={\sf
P}\big(W_1<{\textstyle\frac18}yS_{\gamma,1}^2\big)=1-\sint_{\!\!0}^{\infty}\exp\big\{-{\textstyle\frac18}
yz^2\big\} dG_{\gamma,1}(z),\ \ \ y\ge0.\eqno(14)
$$

The following result generalizes lemma 4 and establishes that the
Weibull distribution with an arbitrary positive shape parameter
$\gamma$ is a scale mixture of the Weibull distribution with an
arbitrary positive shape parameter $\delta>\gamma$.

\smallskip

{\sc Theorem 2}. {\it Let $\delta>\gamma>0$. Then
$$
W_{\gamma}\eqd W_{\delta}\cdot V_{\alpha}^{1/\delta},
$$
where $\alpha=\gamma/\delta\in(0,1)$ and the random variables on the
right-hand side are independent}.

\smallskip

{\sc Proof}. For any $\delta>\gamma>0$, denoting
$\alpha=\gamma/\delta$ (as this is so, we obviously have
$\alpha\in(0,1)$), for any $x\in\mathbb{R}$, from lemma 4 we obtain
$$
{\sf P}(W_{\gamma}>x)=e^{-x^{\gamma}}=e^{-x^{\delta\alpha}}={\sf
P}(W_{\alpha}>x^{\delta})={\sf
P}(W_1>{\textstyle\frac12}S_{\alpha,1}x^{\delta})=
$$
$$
=\sint_{\!\!0}^{\infty}e^{-\frac12
zx^{\delta}}g_{\alpha,1}(z)dz=\sint_{\!\!0}^{\infty}{\sf
P}\big(W_{\delta}>x({\textstyle\frac{1}{2}}z)^{1/\delta}\big)g_{\alpha,1}(z)dz={\sf
P}(W_{\delta}\cdot V_{\alpha}^{1/\delta}>x),
$$
Q. E. D.

\smallskip

{\sc Remark 1}. If $0<\gamma<\delta\le 2$, then the result of
theorem 2 directly follows from theorem 2.3.1 in
\cite{Zolotarev1983} by virtue of the formal coincidence of the
characteristic function of a strictly stable law with the
complementary distribution function of the Weibull law with the
corresponding parameter (see the proof of lemma 3).

\smallskip

The representation of the Weibull distribution with $\gamma\in(0,1)$
in the form of a mixed exponential distribution (see lemma 4) yields
the following by-product result concerning the explicit
representation of the moments of symmetric or one-sided strictly
stable distributions.

\smallskip

{\sc Corollary 2.} {\it Let $S_{\gamma,1}$ be a positive random
variable with the strictly stable distribution with the
characteristic exponent $\gamma\in(0,1)$ and $\theta=1$. For
$\beta\in(0,\gamma)$ we have}
$$
{\sf
E}S_{\gamma,1}^{\beta}=\frac{2^{\beta}\Gamma(1-\beta/\gamma)}{\Gamma(1-\beta)}.
$$

\smallskip

{\sc Proof}. For any $\delta>-\gamma$ we have ${\sf
E}W_{\gamma}^{\delta}=\Gamma(1+\delta/\gamma)$. For such $\delta$
from lemma 3 it follows that ${\sf
E}W_{\gamma}^{\delta}=2^{\delta}{\sf E}W_1^{\delta}{\sf
E}S_{\gamma,1}^{-\delta}$. Setting $\delta=-\beta>-\gamma$ we obtain
the desired assertion.

\smallskip

The statement of lemma 2 can be rewritten as $S_{\alpha,0}\eqd
X\sqrt{S_{\alpha/2,1}}$ with the random variables on the right-hand
side being independent. Therefore, corollary 2, in turn, implies

\smallskip

{\sc Corollary 3.} {\it Let $S_{\alpha,0}$ be a random variable with
the symmetric strictly stable distribution with the characteristic
exponent $\alpha\in(0,2)$ $($see $(6))$. Then for $\beta<\alpha<2$
we have}
$$
{\sf E}|S_{\alpha,0}|^{\beta}=\frac{2^{\beta}}{\sqrt{\pi}}\cdot
\frac{\Gamma\big((\beta+1)/2\big)\Gamma\big((\alpha-\beta)/\alpha\big)}{\Gamma\big((2-\beta)/\beta\big)}.
$$

\section{Product representations for random variables with symmetric
two-sided Weibull distributions}

Let $\gamma>0$. The distribution of the random variable $\widetilde
W_{\gamma}$:
$$
{\sf P}\big(\widetilde W_{\gamma}<x\big)={\textstyle\frac12}
e^{-|x|^{\gamma}}{\bf
1}(x<0)+\big[1-{\textstyle\frac12}e^{-x^{\gamma}}\big]{\bf
1}(x\ge0)\eqno(15)
$$
is called the {\it symmetric two-sided Weibull distribution} with
shape parameter $\gamma$. Distribution (15) was introduced in the
paper \cite{Sornette2000} as a heavy-tailed model for the evaluation
of financial risks. Further generalizations and references can be
found, for example, in the recent works \cite{ChenGerlach2011,
ChenGerlach2013}.

It is easy to see that if $W_{\gamma}$ is a random variable with the
usual (one-sided) Weibull distribution (1) and $U$ is a random
variable taking values $-1$ and $+1$ with probabilities $\frac12$
each and independent of $W_{\gamma}$, then $\widetilde
W_{\gamma}\eqd U W_{\gamma}$ and hence, $|\widetilde W_{\gamma}|\eqd
W_{\gamma}$.

Moreover, from theorem 1 it obviously follows that for
$\gamma\in(0,1]$ we have
$$
\widetilde W_{\gamma}\eqd X \sqrt{2W_1V_{\gamma}^2}
$$
where the random variables on the right-hand side are independent.

\smallskip

{\sc Corollary 4.} {\it For any $\gamma\in(0,1]$ the symmetric
two-sided Weibull distribution with parameter $\gamma$ is a scale
mixture of normal laws}:
$$
{\sf P}(\widetilde W_{\gamma}<x)={\sf
E}\Phi\big(x/\sqrt{2W_1V_{\gamma}^2}\big)=\sint_{\!\!0}^{\infty}\Phi\big(x/\sqrt{y}\big)dH_{\gamma}(y),\
\ \ x\in\mathbb{R},
$$
{\it where, as in $(14)$,}
$$
H_{\gamma}(y)=1-\sint_{\!\!0}^{\infty}\exp\big\{-{\textstyle\frac18}
yz^2\big\} dG_{\gamma,1}(z),\ \ \ y\ge0.
$$

\smallskip

Note that here the mixing distribution $H_{\gamma}(x)$ is exactly
the same as in theorem 1, furthermore, this distribution is
absolutely continuous. It is easy to see that the corresponding
density has the form
$$
h_{\gamma}(x)={\textstyle\frac18}\sint_{\!\!0}^{\infty}z^2\exp\big\{-{\textstyle\frac18}xz^2\big\}g_{\gamma,1}(z)dz,\
\ \ x>0.
$$
Despite the fact that, in general, the mixing distribution
$H_{\gamma}(x)$ cannot be expressed via elementary functions, it is
possible to trace the asymptotic behavior of its tail. Namely, in
\cite{AntonovKoksharov2006} the following result was proved
connecting the tail behavior of a normal scale mixture with that of
the corresponding mixing distribution.

\smallskip

{\sc Lemma 5} \cite{AntonovKoksharov2006}. {\it Assume that a
distribution function $F(x)$ is a scale mixture of zero-mean normal
laws,
$$
F(x)=\sint_{\!\!0}^{\infty}\Phi\Big(\frac{x}{\sqrt{u}}\Big)dQ(y),\ \
\ x\in\R,
$$
with $Q(0)=0$. Let $L(x)$ be a positive slowly varying function,
that is, $L:\mathbb{R}\to\mathbb{R}_+$ so that for any $p>0$
$$
\lim_{x\to\infty}\frac{L(px)}{L(x)}=1.
$$
Let $\rho\in(0,2)$, $\beta>0$. Then
$$
\liminf_{x\to\infty}\frac{\ln[1-F(x)]}{x^{\rho}L(x)}=-\frac{1}{\beta}
$$
if and only if}
$$
\liminf_{u\to\infty}\frac{\ln[1-Q(u)]}{u^{\frac{\rho}{2-\rho}}
\big[L\big(u^{\frac{\rho}{2-\rho}}\big)\big]^{\frac{2}{2-\rho}}}=
-\frac{1}{2\beta^{\frac{2}{2-\rho}}}.
$$

\smallskip

Let $\gamma\in(0,1]$. We obviously have ${\sf P}(\widetilde
W_{\gamma}>x)=\frac12 e^{-x^{\gamma}}$. Therefore from corollary 4
and lemma 5 with $\rho=\gamma$, $\beta=1$ and $L(x)\equiv 1$ it
follows that
$$
1-H_{\gamma}(x)\sim
\exp\big\{-{\textstyle\frac12}x^{\frac{\gamma}{2-\gamma}}\big\},\ \
\ x\to\infty.
$$
that is, the tail of the mixing distribution also decreases as an
exponential power function.

\smallskip

Obviously, $\widetilde W_1$ is the random variable with the Laplace
distribution
$$
L(x)\equiv{\sf P}(\widetilde W_1<x)={\textstyle\frac12} e^x{\bf
1}(x<0)+\big[1-{\textstyle\frac12} e^{-x}\big]{\bf 1}(x\ge0).
$$
It can be made sure that the product representation
$$
\widetilde W_1\eqd X\sqrt{2W_1},
$$
holds, where the random variables on the right-hand side are
independent (see, e. g., \cite{KorolevBeningShorgin2011}, p.
578-579). Then corollary 4, in turn, implies

\smallskip

{\sc Corollary 5.} {\it For any $\gamma\in(0,1]$ the symmetric
two-sided Weibull distribution with parameter $\gamma$ is a scale
mixture of Laplace distributions}:
$$
{\sf P}(\widetilde W_{\gamma}<x)={\sf
E}L\big({\textstyle\frac12}xS_{\gamma,1}\big)=\sint_{\!\!0}^{\infty}L({\textstyle\frac12}xy)g_{\gamma,1}(y)dy,\
\ \ x\in\mathbb{R}.
$$

\smallskip

Actually corollary 5 is a particular case of the following more
general statement which is an analog of theorem 2 for symmetric
two-sided Weibull distributions.

\smallskip

{\sc Theorem 3.} {\it Let $\delta>\gamma>0$. Then
$$
\widetilde W_{\gamma}\eqd \widetilde W_{\delta}\cdot
V_{\alpha}^{1/\delta},
$$
where $\alpha=\gamma/\delta\in(0,1)$, and the random variables on
the right-hand side are independent}.

\smallskip

{\sc Proof.} As it has been noted above, $|\widetilde
W_{\gamma}|\eqd W_{\gamma}$ for any $\gamma>0$. Then by theorem 2
$$
|\widetilde W_{\gamma}|\eqd W_{\gamma}\eqd W_{\delta}\cdot
V_{\alpha}^{1/\delta}\eqd |\widetilde W_{\delta}|\cdot
V_{\alpha}^{1/\delta}= |\widetilde W_{\delta}\cdot
V_{\alpha}^{1/\delta}|.
$$
Hence, since the distributions of $\widetilde W_{\gamma}$ and
$\widetilde W_{\delta}$ are symmetric, for arbitrary $x>0$ we have
$$
{\sf P}(\widetilde W_{\gamma}<x)={\textstyle\frac12}\big[{\sf
P}(|\widetilde W_{\gamma}|<x)+1\big]={\textstyle\frac12}\big[{\sf
P}(|\widetilde W_{\delta}\cdot V_{\alpha}^{1/\delta}|<x)+1\big]={\sf
P}(\widetilde W_{\delta}\cdot V_{\alpha}^{1/\delta}<x).
$$
And if $x<0$, then $-x>0$ so that according to what has already been
proved,
$$
{\sf P}(\widetilde W_{\gamma}<x)=1-{\sf P}(\widetilde
W_{\gamma}<-x)=1-{\sf P}(\widetilde W_{\delta}\cdot
V_{\alpha}^{1/\delta}<-x)={\sf P}(\widetilde W_{\delta}\cdot
V_{\alpha}^{1/\delta}<x).
$$
The theorem is proved.

\section{Formal asymmetric two-sided Weibull distribution as
a scale-location mixture of normal laws}

We begin this section with the following formal definition. Let
$a_1$ and $a_2$ be two finite positive numbers, $\gamma>0$.

\smallskip

{\sc Definition 1}. A random variable $W_{a_1,a_2;\gamma}$ will be
said to have the {\it asymmetric Weibull distribution of the first
kind $\mathfrak{W}_{I}(x)$ with parameters $a_1$, $a_2$, $\gamma$},
if its distribution function has the form
$$
\mathfrak{W}_{I}(x)={\sf P}(W_{a_1,a_2;\gamma}<x) =
\begin{cases}{\displaystyle\frac{a_1 }{a_1+a_2}\cdot e^{-(a_2|x|)^{\gamma}}}, & x\leq0 \vspace{2mm}\\
{\displaystyle 1-\frac{a_2}{a_1+a_2}\cdot e^{-(a_1x)^{\gamma}}}, &
x> 0.
\end{cases}\eqno(16)
$$

\smallskip

Each of positive and negative branches of the distribution
$\mathfrak{W}_{I}(x)$ formally coincides with the classical Weibull
distribution. Our aim in this section is to obtain normal mixture
representations for distribution (16). We will construct this
representation in two stages. On the first stage we will construct
an asymmetric normal-mixture-type representation for the Laplace
(two-sided exponential) distribution.

Let $a_1$ and $a_2$ be two finite positive numbers. We will say that
a random variable $\Lambda_{a_1,a_2}$ has the {\it asymmetric
Laplace distribution with parameters $a_1$ and $a_2$}, if
$\Lambda_{a_1,a_2}\eqd W_{a_1,a_2;1}$, where the distribution of
$W_{a_1,a_2;1}$ is defined by relation (16) with $\gamma=1$.

It is easy to see that the probability density $\ell_{a_1,a_2}(x)$
corresponding to the distribution function $L_{a_1,a_2}(x)={\sf
P}(\Lambda_{a_1,a_2}<x)$ has the form
$$
\ell_{a_1,a_2}(x) = \begin{cases}{\displaystyle\frac{a_1
a_2}{a_1+a_2}\cdot e^{a_2
x}}& \text{for $x\le 0$}\vspace{2mm}\\
\displaystyle{ \frac{a_1 a_2}{a_1+a_2}\cdot e^{-a_1 x}}&\text{for
$x>0$}.\end{cases}
$$
The asymmetric Laplace distribution is a popular model widely used
in many fields, see, e. g., \cite{KotzKozubowskiPodgorski2001}.

Show that this distribution is a special variance-mean normal
mixture.

\smallskip

{\sc Lemma 6.} {\it Let $\mu\in\r$, $\sigma^2\in(0,\infty)$,
$\lambda\in(0,\infty)$. Assume that the random variable $Y$ is
representable in the form
$$
Y\eqd \frac{\sigma}{\sqrt{\lambda}}\cdot
X\sqrt{W_1}+\mu\cdot\frac{W_1}{\lambda},
$$
where the random variable $X$ has the standard normal distribution
and the random variable $W_1$ has the standard exponential
distribution $($so that $W_1/\lambda$ has the exponential
distribution with parameter $\lambda)$, and the random variables $X$
and $W_1$ are independent. Then $Y\eqd\Lambda_{a_1,a_2}$, that is,
$$
{\sf P}(Y<x)={\sf E}\Phi\bigg(\frac{\lambda x-\mu
W_1}{\sigma\sqrt{\lambda W_1}}\bigg)=L_{a_1,a_2}(x),\ \ \ x\in\R,
$$
where}
$$
a_1=\frac{1}{\sqrt{\mu^2+2\lambda\sigma^2}+\mu},\ \ \
a_2=\frac{1}{\sqrt{\mu^2+2\lambda\sigma^2}-\mu}.
$$

\smallskip

{\sc Proof}. It is easy to see that, by the independence of $X$ and
$W_1$, the characteristic function of $Y$ has the form
$$
{\sf E}e^{itY}= {\sf E}e^{it(\sigma X\sqrt{U}+\mu U)}
=\lambda\int_0^{\infty} \exp \big\{z\big( it\mu -
{\textstyle\frac12}\sigma^2t^2 - \lambda\big)\big\} d z =
\frac{\lambda}{it\mu - \frac12 \sigma^2t^2 -\lambda}, \ \ \
t\in\r.\eqno(17)
$$
It remains to make sure that characteristic function (17)
corresponds to the asymmetric Laplace distribution. To obtain a
simpler form of the right-hand side of (17), change the parameters:
$$
\begin{cases}{\displaystyle w-v = \frac{\mu }{\lambda}}, & \vspace{1mm}\\
{\displaystyle v\cdot w  = \frac{\sigma^2}{2\lambda}}. &
\end{cases}\eqno(18)
$$
From the first equation (18) we obtain
$$
w=v+\frac{\mu }{\lambda},
$$
whereas the second equation yields
$$
v\Big(v+\frac{\mu }{\lambda}\Big) = \frac{\sigma^2}{2\lambda}
$$
or
$$
v^2 + \frac{\mu }{\lambda}v - \frac{\sigma^2}{2\lambda} = 0.
$$
System (18) has two solutions with respect to $v$:
$$
v_1 = -\frac{\mu }{2\lambda} + \frac12\sqrt{\frac{\mu
^2}{\lambda^2}+\frac{2\sigma^2}{\lambda}},\ \ \ \ v_2 = -\frac{\mu
}{2\lambda} - \frac12\sqrt{\frac{\mu
^2}{\lambda^2}+\frac{2\sigma^2}{\lambda}}.
$$
One of them, $v_1$, is positive. As this is so, $w_1 = v_1
+\frac{\mu }{\lambda}$ is also positive. The parameters $v=v_1$ and
$w=w_1$ will be used in what follows.

With the new parametrization, characteristic function (17) takes the
form
$$
{\sf E}e^{itY}=\frac{1}{(1-iwt)}\cdot \frac{1}{(1+ivt)}.
$$
Notice that $\frac{1}{1-iwt}$ is the characteristic function of the
exponential distribution with parameter $a_1 = \frac1w$
corresponding to the density
$$
p_1(x) = \begin{cases} 0, & x\leq 0\vspace{1mm}\\
a_1 e^{-a_1 x}, & x>0.
\end{cases}
$$
At the same time $\frac{1}{1+ivt}$ is the characteristic function
corresponding to the density
$$
p_2(x) = \begin{cases}a_2 e^{a_2 x}, & x\leq0 \vspace{1mm}\\
0, & x> 0,
\end{cases}
$$
where $a_2 = \frac1v.$ Hence, $g(t)$ is the characteristic function
of the convolution $p(x)$ of the densities $p_1(x)$ и $p_2(x)$ which
has the following form: for $x\leq0$
$$
p(x) = \int_{-\infty}^{x}a_1 e^{-a_1(x-y)}a_2 e^{a_2 y} dy =
\frac{a_1 a_2}{a_1+a_2}e^{a_2 x},
$$
and for $x>0$
$$
p(x) = \int_{-\infty}^{0}a_1 e^{-a_1(x-y)}a_2 e^{a_2 y} dy =
\frac{a_1 a_2}{a_1+a_2}e^{-a_1 x},
$$
that is, $p(x)=\ell_{a_1,a_2}(x)$, $x\in\R$. In other words, this
density corresponds to asymmetric Laplace distribution (1).
Returning to the original parameters we obtain
$$
a_1=\big(\sqrt{\mu ^2+2\lambda\sigma^2}+\mu \big)^{-1},\ \ \
a_2=\big(\sqrt{\mu ^2+2\lambda\sigma^2}-\mu \big)^{-1}.
$$
The lemma is proved.

\smallskip

From lemma 4 by formal calculation it follows that for
$\gamma\in(0,1]$ and $y\ge 0$
$$
{\sf P}(\Lambda_{a_1,a_2}\cdot
V_{\gamma}>y)=\frac{a_2}{a_1+a_2}\sint_{\!\!0}^{\infty}e^{-a_1yz}g_{\gamma,1}(z)dz=\frac{a_2e^{-(a_1y)^{\gamma}}}{a_1+a_2},\eqno(19)
$$
and for $\gamma\in(0,1]$ and $y<0$
$$
{\sf P}(\Lambda_{a_1,a_2}\cdot
V_{\gamma}<y)=\frac{a_1}{a_1+a_2}\sint_{\!\!0}^{\infty}e^{-a_2|y|z}g_{\gamma,1}(z)dz=\frac{a_1e^{-(a_2|y|)^{\gamma}}}{a_1+a_2}.\eqno(20)
$$

Now from lemma 5 and relations (19) and (20) we obtain the following
statement.

\smallskip

{\sc Theorem 3.} {\it Let $\gamma\in(0,1]$, $\mu\in\r$,
$\sigma^2\in(0,\infty)$, $\lambda\in(0,\infty)$. Assume that a
random variable $Z$ is representable in the form
$$
Z\eqd \Big(\frac{\sigma}{\sqrt{\lambda}}\cdot X\sqrt{W_1}+\frac{\mu
W_1}{\lambda}\Big)\cdot V_{\gamma},
$$
where the random variable $X$ has the standard normal distribution,
the random variable $W_1$ has the standard exponential distribution
$($so that $W_1/\lambda$ has the exponential distribution with
parameter $\lambda)$, $V_{\gamma}=2S_{\gamma,1}^{-1}$, where
$S_{\gamma,1}$ is the random variable with the one-sided strictly
stable density $g_{\gamma,1}(x)$, moreover, $X$, $W_1$ и
$S_{\gamma,1}$ are independent. Then $Z\eqd W_{a_1,a_2;\gamma}$,
that is,
$$
{\sf P}(Z<x)={\sf E}\Phi\bigg(\frac{\lambda x-\mu W_1
V_{\gamma}}{\sigma\sqrt{\lambda W_1}V_{\gamma}}\bigg)={\sf
P}(W_{a_1,a_2;\gamma}<x)=\mathfrak{W}_{I}(x),\ \ \ x\in\R,\eqno(21)
$$
where}
$$
a_1=\frac{1}{\sqrt{\mu^2+2\lambda\sigma^2}+\mu},\ \ \
a_2=\frac{1}{\sqrt{\mu^2+2\lambda\sigma^2}-\mu}.
$$

\section{Asymmetric generalization of the two-sided Weibull distribution
by variance-mean mixing of normal laws}

Although the asymmetric generalization $\mathfrak{W}_{I}(x)$ of the
two-sided Weibull distribution introduced in the preceding section
is a scale-location mixture of normal laws, it is not so easy to
give an example of a simple limit scheme, say, for random sums of
independent random variables with such a distribution as a limit
law, since the random shift and scale parameters in (21) are linked
in a non-trivial way. However, if we change the definition of the
asymmetric Weibull distribution in a reasonable way, then such a
limit scheme can be constructed rather easily. For this purpose we
will use a representation of the Weibull-type distribution as a
variance-mean normal mixture based on the representation of the
symmetric two-sided Weibull law obtained in corollary 4.

\smallskip

{\sc Definition 2.} Let $\mu\in\R$, $\sigma>0$, $\gamma\in(0,1]$. A
random variable $W^*_{\mu,\sigma;\gamma}$ will be said to have the
{\it asymmetric Weibull distribution of the second kind
$\mathfrak{W}_{II}(x)$ with parameters $\mu$, $\sigma$ and
$\gamma$}, if its distribution function has the form of the
variance-mean normal mixture
$$
\mathfrak{W}_{II}(x)=\sint_{\!\!0}^{\infty}\Phi\Big(\frac{x-\mu
z}{\sigma\sqrt{z}}\Big)dH_{\gamma}(z),\ \ \ x\in\R,
$$
where the mixing distribution function $H_{\gamma}(z)$ is defined in
(14).

\smallskip

From corollary 4 it obviously follows that with $\mu=0$, the
asymmetric two-sided Weibull distribution of the second kind
$\mathfrak{W}_{II}(x)$ has the form (15). But if $\mu\neq0$, then
this distribution cannot be expressed in terms of elementary
functions. Nevertheless, it is rather easy to formulate a limit
theorem for random sums of independent identically distributed
random variables with finite variances in which the asymmetric
two-sided Weibull distribution of the second kind
$\mathfrak{W}_{II}(x)$ turns out to be the limit law.

Let $\{X_{n,j}\}_{j\ge1}$, $n=1,2,\ldots$, be double array of
row-wise independent and identically distributed random variables.
Let $\{N_n\}_{n\ge1}$ be a sequence of nonnegative integer-valued
random variables such that for each $n\ge1$ the random variables
$N_n,X_{n,1},X_{n,2},\ldots$ are independent. For any
$n,k\in\mathbb{N}$ let
$$
S_{n,k}=X_{n,1}+\ldots +X_{n,k}.
$$
To avoid misunderstanding, assume that $\sum_{j=1}^0=0$. The symbol
$\Longrightarrow$ will denote the convergence in distribution.

In \cite{Korolev2013} the following statement was proved.

\smallskip

{\sc Lemma 7}. {\it Assume that there exist a sequence of natural
numbers $\{k_n\}_{n\ge1}$ and numbers $\mu\in\mathbb{R}$ and
$\sigma>0$ such that
$$
{\sf P}\big(S_{n,k_n}<x\big)\Longrightarrow
\Phi\Big(\frac{x-\mu}{\sigma}\Big).\eqno(22)
$$
Assume that $N_n\to\infty$ in probability. Then the distributions of
the random sums $S_{N_n}$ converge to some distribution function
$F(x):$
$$
{\sf P}\big(S_{n,N_n}<x\big)\Longrightarrow F(x),
$$
if and only if there exists a distribution function $Q(x)$ such that
$Q(0)=0$,
$$
F(x)=\sint_{\!\!0}^{\infty}\Phi\Big(\frac{x-\mu
z}{\sigma\sqrt{z}}\Big)dQ(z),
$$
and
$$
{\sf P}(N_n<xk_n)\Longrightarrow Q(x).
$$
}

\smallskip

{\sc Remark 2.} Condition (22) holds in the following rather general
situation. Assume that the random variables $X_{n,j}$ have finite
variances. Also assume that for each $n$ and $j$
$$
X_{n,j}=X_{n,j}^*+\mu_n,
$$
where $\mu_n\in\mathbb{R}$ and $X_{n,j}^*$ is a random variable with
${\sf E} X_{n,j}^*=0$, ${\sf D} X_{n,j}^*=\sigma_n^2<\infty$, so
that ${\sf E} X_{n,1}=\mu_n$ и ${\sf D} X_{n,1}=\sigma_n^2$. Assume
that $\mu_nk_n\to \mu\in\R$ и $k_n\sigma_n^2\to
\sigma^2\in(0,\infty)$ as $n\to\infty$. Then according to the
classical result on necessary and sufficient conditions for the
convergence of the distributions of independent identically
distributed random variables with finite variances to the normal law
in the double array limit scheme (see, e. g.,
\cite{GnedenkoKolmogorov1949}), we can see that convergence (22)
takes place if and only if the Lindeberg condition holds:
$$
\lim_{n\to\infty}k_n{\sf
E}(X_{n,1}^*)^2\mathbf{1}(|X_{n,1}^*|\ge\varepsilon)=0
$$
for any $\varepsilon>0$.

\smallskip

From lemma 7 and definition 2 we immediately obtain the following
statement establishing necessary and sufficient conditions for the
convergence of the distributions of random sums of independent
identically distributed random variables to the asymmetric two-sided
Weibull distribution of the second kind $\mathfrak{W}_{II}(x)$.

\smallskip

{\sc Theorem 4.} {\it Assume that there exist a sequence of natural
numbers $\{k_n\}_{n\ge1}$ and numbers $\mu\in\mathbb{R}$ and
$\sigma>0$ such that convergence $(22)$ takes place. Assume that
$N_n\to\infty$ in probability. Then the distributions of random sums
$S_{N_n}$ of independent identically distributed random variables
converge to to the asymmetric two-sided Weibull distribution of the
second kind $\mathfrak{W}_{II}(x)$ with parameters $\mu$, $\sigma$,
$\gamma:$
$$
{\sf P}\big(S_{n,N_n}<x\big)\Longrightarrow \mathfrak{W}_{II}(x),
$$
if and only if
$$
{\sf P}(N_n<xk_n)\Longrightarrow H_{\gamma}(x),
$$
where the distribution function $H_{\gamma}(x)$ is defined in} (14).

\renewcommand{\refname}{References}

\end{document}